\documentclass[a4paper, reqno]{amsart}
\usepackage{a4wide}

\usepackage[english]{babel}
\usepackage[utf8]{inputenc}

\usepackage{amsmath}
\usepackage{amsthm}
\usepackage{amssymb} 

\usepackage{mathtools}
    \mathtoolsset{centercolon}

\usepackage{xcolor}

\usepackage{enumitem}
\usepackage{hyperref}
\usepackage[nameinlink,capitalize]{cleveref}

\usepackage[alphabetic,backrefs]{amsrefs}

\usepackage{tikz}

\theoremstyle{plain}
  
  \newtheorem{thm}{Theorem}[section]
  
  \newtheorem{prop}[thm]{Proposition}
  \newtheorem{lemma}[thm]{Lemma}

\theoremstyle{definition}
  \newtheorem{defn}[thm]{Definition}
  \newtheorem{nota}[thm]{Notation}
  \newtheorem{convention}[thm]{Convention}
  \newtheorem{remark}[thm]{Remark}
  \newtheorem{example}[thm]{Example}
  \newtheorem{question}[thm]{Question}

\DeclareMathOperator{\Aut}{Aut}

\DeclareMathOperator{\id}{id}
\DeclareMathOperator{\Inn}{Inn}
\DeclareMathOperator{\nucl}{nucl}
\DeclareMathOperator{\GL}{GL}
\DeclareMathOperator{\SL}{SL}
\DeclareMathOperator{\SO}{SO}
\DeclareMathOperator{\SU}{SU}

\newcommand{\CC}{\mathbb{C}}
\newcommand{\FF}{\mathbb{F}}
\newcommand{\NN}{\mathbb{N}}
\newcommand{\RR}{\mathbb{R}}
\newcommand{\ZZ}{\mathbb{Z}}

\newcommand{\PPP}{\mathcal{P}}

\let\ol\overline

\newcommand\eps\varepsilon

\DeclarePairedDelimiter\gen{\langle}{\rangle}

\newcommand{\mtr}[4]{\begin{pmatrix} {#1} & {#2} \\ {#3} & {#4} \end{pmatrix}} 

\newcommand{\Defn}[1]{\textcolor{blue}{\textit{#1}}}

\setenumerate[1]{label={(\roman*)}}

\begin{document}

\title{Decompositions of Kac-Moody groups}
\author{Max Horn}
\maketitle

\section{Introduction}

Let $G$ be a split (minimal) Kac-Moody group over $\RR$ or $\CC$ with maximal torus $T$, and let $\theta$ be a Cartan-Chevalley involution of $G$, twisted by complex conjugation, and satisfying that $\theta(T)=T$. Furthermore, let $K$ be the subgroup fixed by $\theta$, and $\tau:G\to G, g\mapsto g\theta(g)^{-1}$. Let $A:=\tau(T)$.
 
In this note, we show resp. revisit that $G$ admits a (refined) Iwasawa decompositions $G=UAK$. We also show that if $G$ is of non-spherical type, then it never admits a polar decomposition $G=\tau(G)K$ nor a Cartan decompositions $G=KAK$. This has implications for the geometrical structure of the Kac-Moody symmetric space $G/K \cong \tau(G)$ as defined and studied in \cite{FHHK}.

\medskip

\textbf{Acknowledgements.}
I would like to thank Walter Freyn, Tobias Hartnick and Ralf Köhl for many inspiring discussions on Kac-Moody symmetric spaces, motivating me to write this note.

\section{Basics}

Throughout this note, we assume that the reader is familiar with topics such as Kac-Moody groups, twin buildings, and so on. A brief summary of the required theory, close in notation to what we use here, can be found in \cite{FHHK}*{Section 3} (see also \cite{Gramlich/Horn/Muehlherr}). For a comprehensive reference, we refer to \cite{Abramenko/Brown:2008}.

\begin{nota}
Throughout this paper, let $G$ be a split (minimal) Kac-Moody group of rank $n$ over some field $\FF$. We fix the following notation:
\begin{itemize}
\item $(W,S)$: the associated Coxeter system, with $W$ the Weyl group of $G$.

\item $\Phi$ is the associated root system, with $\Pi=\{\alpha_1,\dots,\alpha_n\}$ a system of fundamental roots,
and corresponding sets $\Phi_+$ resp. $\Phi_-$ of positive resp. negative roots.

\item $\{U_\alpha\}_{\alpha\in\Phi}$ is a root group datum for $G$ (cf. \cite{CR09}).

\item $T:=\cap_{\alpha\in\Phi} N_G(U_\alpha)$ is a maximal torus of $G$.

\item $U_\eps:=\gen{U_\alpha \mid \alpha \in \Phi_\eps}$ for $\eps\in\{+,-\}$.

\item $B_\eps:=TU_\eps$ for $\eps\in\{+,-\}$ is a Borel subgroup of $G$, with unipotent radical $U_\eps$.

\item $(B_+,B_-,N)$: the associated twin $BN$-pair.
\item $G_\alpha:=\gen{ U_\alpha, U_{-\alpha} }$ for $\alpha\in\Pi$ is a fundamental rank-1-subgroup; since $G$ is split, $G_\alpha$ is isomorphic to a central quotient of $\SL_2(\FF)$.

\item $\Delta:=((\Delta_+,\delta_+),(\Delta_-,\delta_-),\delta_*)$ is the twin building associated to $G$ (we identify $\Delta_\eps$, when viewed as a chamber system, with $G/B_\eps$).

\end{itemize}
\end{nota}

\begin{remark}
\begin{enumerate}
\item One has $U_+\cap U_- =\{1\}$, and $B_+\cap B_- = T$.

\item $G$ is generated by the root groups $U_\alpha$ and the torus $T$.
\end{enumerate}
\end{remark}

\begin{example} \label{example KM groups}
\begin{enumerate}
\item Let $n\geq 1$ and $G=\SL_{n+1}(\FF)$. This is a split Kac-Moody group of type $A_n$. Here $B_\eps$ are the subgroups of upper resp. lower triangular matrices; $T$ the subgroup of diagonal matrices; $U_\eps$ the subgroups of strictly upper resp. lower triangular matrices.
The Weyl group then is isomorphic to $S_n$, and of type $A_n$, as is therefore each half of the twin building.

More generally, any split reductive algebraic group is an example.

\item 
However, we are mainly interested in the \Defn{non-spherical} case, that is, when $W$ is infinite. 
As an example for this, consider $G=\SL_n(\FF[t,t^{-1}])$, for some $n\geq 2$.
This is of type $\tilde{A}_{n-1}$ and rank $n$.

\end{enumerate}
\end{example}

\begin{defn}[See \cite{Caprace:2009}]
Let $g\in G$.
\begin{enumerate}
\item
We call $g$ \Defn{diagonalizable} if it conjugate to an element of $T$. Equivalently, it stabilizes a pair of opposite chambers in the twin building $\Delta$, and hence stabilizes the twin apartment spanned by them.
\item We call $g$ \Defn{bounded} if it stabilizes spherical residues in each half of the twin building associated to $G$.
\end{enumerate}
\end{defn}

\begin{defn}
Let $\sigma\in\Aut(\FF)$ with $\sigma^2=\id$. A \Defn{$\sigma$-twisted Cartan-Chevalley-involution} of $G$ is an automorphism of $G$ which is $\Inn(G)$-conjugate to an involution $\theta\in\Aut(G)$ satisfying the following for all $\alpha\in\Phi$:
\begin{enumerate}
\item $U_\alpha^\theta=U_{-\alpha}$,
\item $\theta\circ\sigma$ induces a Cartan-Chevalley involution on $G_\alpha$.
\end{enumerate}
\end{defn}

\begin{remark}
Let $\theta$ be a twisted $\sigma$-twisted Cartan-Chevalley-involution.
\begin{enumerate}
\item
Since conjugation by $G$ resp. $\Inn(G)$ changes nothing for the results of interest for us, we will from now on simply assume that $\theta$ has the properties (i) and (ii).

\item
By \cite{Medts/Gramlich/Horn:2009}*{Lemma 4.2} and the discussing preceding it,
restricting $\theta$ to $G_\alpha$ yields an automorphism induced by a map of the form
\[ \SL_2(\FF)\to\SL_2(\FF),\ x\mapsto 
\mtr{0}{1}{-\eps}{0} x^\sigma \mtr{0}{-\eps^{-1}}{1}{0}
= \mtr100{\eps} ((x^\sigma)^T)^{-1} \mtr100{\eps^{-1}}
\]
where $\eps\in\FF$ satisfies $\eps^\sigma=\eps$.
In other words, $\theta$ locally splits into a field automorphism, a Cartan-Chevalley automorphism (also known as sign automorphism), and a diagonal automorphism.

\item
In \cite{Caprace:2009} (see also \cite{Caprace/Muehlherr:2005,Caprace/Muehlherr:2006}), the isomorphism problem of Kac-Moody groups is solved. It turns out that every automorphism of $G$ is the product of an inner automorphism, a field automorphism, a sign automorphism, a diagonal automorphism, and a graph automorphism.

The definition of Cartan-Chevalley-involution, plus our preceding assumption, implies that $\theta$ globally splits into the product of a field automorphism, a sign automorphism, and a diagonal automorphism.

\end{enumerate}
\end{remark}

\begin{nota} From now on, $\theta$ will be a $\sigma$-twisted Cartan-Chevalley-involution, and we fix the following notation:
\begin{itemize}
 \item $K:=\{ g\in G \mid \theta(g)=g \}$, the \Defn{unitary form} of $G$.
 \item $Q:=\{ g\in G \mid \theta(g)=g^{-1} \}$ is the set of \Defn{$\theta$-symmetric} elements in $G$.
 \item $\tau:G\to Q: g\mapsto g\theta(g)^{-1}$ is the \Defn{twist map} associated to $\theta$.
 \item $A:= \tau(T)$ is a \Defn{maximal flat}.
 \item $M:= K\cap T$.
\end{itemize}
\end{nota}

\begin{remark}
\begin{enumerate}
\item $\theta$ induces the inversion map $g\mapsto g^{-1}$ on $Q$ and on $\tau(G) \subseteq Q$.
\item We have $\theta(B_+)=B_-$ and $\theta(B_-)=B_+$, hence $\theta(T)=T$. 

\item $\theta$ induces an involutory bijection between $\Delta_+=G/B_+$ and $\Delta_-=G/B_-$ via $gB_+\mapsto\theta(gB_+)=\theta(g)B_-$. We will refer to this map also as $\theta$. Note that $\theta$ preserves the Weyl (co-)distances $\delta_+$, $\delta_-$ and $\delta_*$. See also \cite{Gramlich/Horn/Muehlherr}*{Proposition 3.1}.

\end{enumerate}
\end{remark}

\begin{lemma}
The restriction of $\theta$ to $T$ is the map $t\mapsto \sigma(t)^{-1}$. Hence $A=\{ tt^\sigma \mid t\in T\}$.
\end{lemma}

\begin{proof}
This follows from the fact $\theta$ decomposes into the product of a diagonal automorphism, a field automorphism and  sign automorphism. The diagonal automorphism acts trivially on $T$ and the sign automorphism acts by inversion.
\end{proof}

\begin{example} \label{example CC invs}
We continue Example~\ref{example KM groups}.
\begin{enumerate}
\item Let $G=\SL_{n+1}(\RR)$, and consider the Cartan-Chevalley involution $\theta(g):=(g^T)^{-1}$ on $G$.
Then $K$ is the special orthogonal group $\SO_{n+1}(\RR)$, and $Q$ the set of symmetric matrices with determinant 1, and $\tau(g)=gg^T$. Thus, $\tau(G)$ consists of the symmetric positive definite matrices with determinant 1. Finally, $A$ consists of the positive diagonal matrices with determinant 1.

For $\FF=\CC$ and $\sigma$ complex conjugation, a $\sigma$-twisted Cartan-Chevalley involution on $G$ is given by $\theta(g):=((g^\sigma)^T)^{-1}$.
We then have $K=\SU_{n+1}(\RR)$, $Q$ is the set of Hermitian matrices with determinant 1, and $\tau(G)\subset Q$ the subset of positive definite matrices.
$A$ is the set of positive diagonal matrices, and $M$ the set of diagonal matrices with all entries $\pm1$.

\item
Let $G=\SL_n(\FF[t,t^{-1}])$ and $\sigma\in\Aut(\FF)$ with $\sigma^2=\id$.
Then a $\sigma$-twisted Cartan-Chevalley involution on $G$ is given by $\theta(x):=((x^{\sigma\rho})^{-1})^T$, where $\rho$ is the unique $\FF$-linear ring automorphism of $\FF[t,t^{-1}]$ interchanging $t$ and $t^{-1}$.

Here, $A$ is again the set of positive diagonal matrices, and $M$ the set of diagonal matrices with all entries $\pm1$.

\end{enumerate}
\end{example}

In this note, we study the (non-)existence of various decompositions of $G$:

\begin{defn} \label{def decomps}
 $G$ admits, with respect to $\theta$, \ldots
\begin{itemize}
\item \ldots an \Defn{Iwasawa decomposition} if $G=BK$ holds.

\item \ldots a \Defn{refined Iwasawa decomposition} if $K\times A\times U\to G, (k,a,u)\mapsto kau$ is a bijection.

\item \ldots a \Defn{polar decomposition} if $G=\tau(G) K$ holds.

\item \ldots a \Defn{Cartan decomposition} if $G=KAK$ holds.

\item \ldots a \Defn{Kostant decomposition} if $G=KUK$ holds.
\end{itemize}
\end{defn}

\section{Iwasawa decomposition}

\begin{convention}
From now on for the rest of this paper, we will assume that either $\FF=\RR$ and $\sigma=\id$, or else $\FF=\CC$ and $\sigma$ is complex conjugation.
\end{convention}

The existence of a refined Iwasawa decomposition for complex Kac-Moody groups has been known for quite some time, see e.g. \cite{Kac/Peterson:1983}*{Corollary 4}. However, no proofs are given there. The real case is at the very least known as folklore, though I am not aware of a fully developed proof in the literature.

The existence of non-refined Iwasawa decompositions $G=BK$ of $G$ over arbitrary fields was studied extensively in \cite{Medts/Gramlich/Horn:2009}. This actually allows generalizing various results in later sections of this note beyond the real and complex case. Despite this, we mostly focus on the real and complex, as this allows for a particularly simple exposition, and is the case we are currently most interested in for applications, see \cite{FHHK}.

Thus focusing again on the real and complex case, we can rephrase the existence of an Iwasawa decomposition $G=BK$ as saying that the map $B\times K\to G,\ (b,k)\mapsto G$ is surjective. In general, this map is not injective. To rectify this, one may replace $B$ with a suitable subgroup, and study when $G$ admits a refined Iwasawa decomposition as defined above. The existence of such a refined Iwasawa decomposition in the real case is also shown in \cite{FHHK}*{Theorem 3.23}. Virtually the same argument applies for the complex case. For the convenience of the reader, we give a full proof. Note that loc.cit.\ also describes and proves a \emph{topological} Iwasawa decomposition, but only in the real case; the complex case is currently open. 

\begin{lemma}
$M\cap A=\{1\}$ and $T=MA$ hold.
This induces an isomorphism of topological groups $T\cong M\times A$.
\end{lemma}

\begin{proof}
We have $T\cong\FF^*)^m$ for some natural number $m$. As stated before, $\theta$ induces on $T$ inversion, composed with complex conjugation.
If $x\in M\cap A$, then on the one hand, $x\in A = \tau(T)$, so $x = t\overline{t} \in \RR_{>0}^m$.
On the other hand, $x\in K$, hence $x=\theta(x)=\overline{x}^{-1}$, so $x\overline{x}=1$. Together this implies $x=1$.
The claim now follows from the polar decomposition $\CC^* \cong \RR_{>0} \times \{ \rho\in\CC \mid |\rho|=1\} \cong \RR_{>0} \times S^1$ in the complex case, and from $\RR^* \cong \RR_{>0} \times \{\pm 1\} \cong \RR_{>0} \times S^0$ in the real case.
\end{proof}

\begin{lemma}
$K\cap B_\eps=M$ holds for $\eps\in\{+,-\}$.
\end{lemma}

\begin{proof}
Clearly $M\subseteq K\cap B_\eps$. Let $k\in K\cap B_{\eps}$. Then $k=\theta(k)\in K\cap B_{-\eps}$,
hence $k\in K\cap B_+\cap B_- = K\cap T = M$.
\end{proof}

\begin{prop}[Refined Iwasawa decomposition]
For $\eps\in\{+,-\}$, the maps
\[ \mu_\eps: K\times A\times U_\eps \to G,\ (k,a,u) \to kau \]
are bijections.
\end{prop}

\begin{proof}
Surjectivity follows from $M\subseteq K$, the preceding lemmas and the (unrefined) Iwasawa decomposition:
\[ KAU_\eps = KMAU_\eps = KTU_\eps = KB_+ = G. \]
Suppose now $kau = k'a'u'$, then $K\ni (k')^{-1}k = a'u'u^{-1} a^{-1} \in AU_\eps A = B_\eps$, whence  $a'u'u^{-1} a^{-1} \in B_\eps\cap K = M$. Therefore $u'u^{-1} \in U_\eps \cap T = \{1\}$, so $u'=u$. This implies $a'a^{-1}\in M\cap A=\{1\}$, and so $a'=a$. We finally conclude from this $k=k'$, thus $\mu_\eps$ is indeed injective.
\end{proof}

\section{Symmetric elements of $G$}

The existence of Iwasawa decompositions $G=KB_\eps$ for $\eps\in\{+,-\}$ implies that all Borel subgroups, i.e., the $G$-conjugates of $B_+$ and $B_-$, are in fact $K$-conjugate to $B_+$ or $B_-$. In particular, all Borel subgroups $B$ of $G$ are \Defn{$\theta$-split}, i.e., for all $g\in G$ we have $B^g\cap\theta(B^g)$ is maximal torus, conjugate to $T$.
This is one of the many ingredient of the following useful lemma.

\begin{remark}
In view of the example $G=\SL_{n+1}(\RR)$, $\theta(g):=(g^T)^{-1}$, we may think of lemma as a generalization of the observation that every real symmetric matrix can be diagonalized by an orthogonal matrix, resp. every Hermitian matrix can be diagonalized by a unitary matrix. Indeed, we use this explicitly
\end{remark}

\begin{lemma} \label{lem:sym-ss}
If $g\in G$ is $\theta$-symmetric, i.e., if $\theta(g)=g^{-1}$,
the following are equivalent:
\begin{enumerate}
\item\label{enum:theta-apt} $g$ fixes a $\theta$-stable twin apartment chamberwise.
\item\label{enum:apt} $g$ fixes a twin  apartment chamberwise (i.e., is diagonalizable).
\item\label{enum:cham} $g$ stabilizes a chamber.
\item\label{enum:orb-all} For all chambers $d$, the $G$-orbit $\{ g^n.d \mid n\in \ZZ\}$ is bounded in the gallery metric.
\item\label{enum:orb} 
For some chamber $d$, the $G$-orbit $\{ g^n.d \mid n\in \ZZ\}$ is bounded in the gallery metric.

\item\label{enum:res} $g$ stabilizes a spherical residue in either half of the twin building.
\end{enumerate}
\end{lemma}

\begin{proof}
The implications $\ref{enum:theta-apt}\implies\ref{enum:apt}\implies\ref{enum:cham}$ and $\ref{enum:orb-all}\implies\ref{enum:orb}$ are elementary.
\begin{description}
\item[$\ref{enum:cham}\implies\ref{enum:theta-apt}$]
Let $c$ be a chamber stabilized by $g$. Since $\theta(c)$
is opposite $c$ and $\theta(c)= \theta(g.c) = g^{-1} . \theta(c)$,
we conclude that $g$ stabilizes the $\theta$-stable twin apartment
$\Sigma(c,\theta(c))$.

\item[$\ref{enum:cham}\implies\ref{enum:orb-all}$]
Suppose $c\in\Delta_+$ is a chamber stabilized by $g$. Then for $d\in\Delta_+$ and $n\in\ZZ$, we have $\delta_+(c,d)=\delta_+(g^n.c,g^n.d)=\delta_+(c,g^n.d)$. Let $\ell:W\to\NN$ be the length map of $(W,S)$.
Then by the triangle inequality for the building $W$-metric, we have $\ell(\delta_+(d,g^n.d)) \leq \ell(\delta_+(d,c))+\ell(\delta_+(c,g^n.d)) = 2\ell(\delta_+(d,c))$.

\item[$\ref{enum:orb}\implies\ref{enum:res}$] This follows from the
Bruhat-Tits fixed point theorem applied to the CAT(0)-realization of the building; see e.g. \cite{Abramenko/Brown:2008}*{Corollary 12.67}.

\item[$\ref{enum:res}\implies\ref{enum:apt}$]
Let $R$ be a spherical residue stabilized by $g$.
Then we have $\theta(R)=\theta(g.R)=\theta(g).\theta(R) = g^{-1}.\theta(R)$,
and thus $\theta(R)$ is also stabilized by $g$.
In the terminology of \cite{Caprace:2009} this means that $g$ is bounded.

We now consider the automorphism of the spherical building $R$ induced by $g$.
The automorphism group of $R$ is a reductive algebraic group, and can be considered as a subgroup of $\GL_{n+1}(\RR)$. By \cite[Proposition~16.1.5]{HilgertNeeb12}, we can then model $\theta$ as transpose-inverse, composed with complex conjugation (if $\FF=\CC$).
Hence $\theta(g)=g^{-1}$ implies $g^T=\ol{g}$, i.e., $g$ is Hermitian and therefore diagonalizable. Thus it fixes a chamber in $R$, hence in $\Delta_+$.
\qedhere
\end{description}
\end{proof}

\section{The nucleus of $\theta$}

In the next section, we will show that Kac-Moody groups of non-spherical types admit no Polar decomposition $G=\tau(G)K$. To facilitate this, we first collect some general observations about Polar decompositions.

The following two elementary lemmas hold for any group $G$ and involution $\theta\in\Aut(G)$.
\begin{lemma} \label{tauX=tauY iff XK=YK}
 Let $X,Y\subseteq G$. Then $\tau(X)=\tau(Y)$
 if and only if $XK=YK$.
\end{lemma}
\begin{proof}
 For $g,h\in G$, we have
\begin{align*}
   \tau(g)=\tau(h)
 &\iff  g\theta(g)^{-1} = h\theta(h)^{-1} \\
 &\iff  h^{-1}g=\theta(h)^{-1}\theta(g)=\theta(h^{-1}g) \\
 &\iff  h^{-1}g \in K \\
 &\iff  gK = hK. \qedhere
\end{align*}
\end{proof}

\begin{lemma} \label{G=tauK iff tautau=tau}
 $G=\tau(G)K$ holds if and only if $\tau(\tau(G))=\tau(G)$.
\end{lemma}
\begin{proof}
Follows from \cref{tauX=tauY iff XK=YK} for $X:=\tau(G)$, $Y:=G$ and using that $G=GK$.
\end{proof}
Recall that $Q:=\{g\in G \mid \theta(g)=g^{-1} \}$ and $\tau(G)\subseteq Q$.
Then for any $g\in Q$, we have $\tau(g)=g\theta(g^{-1})=g^2$. Thus, $\tau$
acts like the square map on $Q$ and also on $\tau(G)$.
Hence $\tau(G)=\tau(\tau(G))$ is equivalent to requiring that every element of $g\in\tau(G)$
admits a ``square root'', i.e., there is an element $h\in\tau(G)$ such that $\tau(h)=h^2=g$.
But this then also implies that every element has a fourth root, an eighth root and so on.
This motivates the following definition and the subsequent reformulation of the lemma.

\begin{defn}
The \Defn{nucleus} of $X\subseteq G$ is defined as
 \[ \nucl(X) := \bigcap_{k=0}^\infty \left\{x^{(2^k)} \mid x\in X \right\}
 = \{ x\in X \mid \forall k\in\NN\, \exists y\in X : x=y^{2^k} \}
 . \]
\end{defn}

\begin{lemma} \label{G=tau(G)K iff nucl(tau(G))=tau(G)}
$G=\tau(G)K$ holds if and only if $\nucl(\tau(G))=\tau(G)$.
\end{lemma}

\begin{proof}
By \cref{G=tauK iff tautau=tau},  $G=\tau(G)K$ is equivalent to $\tau^2(G)=\tau(G)$,
which in turn is equivalent to $\tau^{n+1}(G)=\tau(G)$ holding for all $n\in\NN$.
Since $\tau^2(g)=\tau(g)\theta(\tau(g))^{-1}=\tau(g)^2$ for all $g\in G$, it follows that $\tau^{n+1}(g)=\tau(g)^{2^n}$.

Hence $\tau^2(G)=\tau(G)$ implies $\tau(G) = \tau^n(G) = \nucl(\tau(G))$ as claimed.
The converse implication follows then from $\tau^2(G)\subseteq \tau(G) = \nucl(\tau(G)) \subseteq \tau^2(G)$.
\end{proof}

\begin{prop}
For $X\subseteq G$, the elements of $\nucl(X)$ are bounded.
\end{prop}

\begin{proof}
Clearly $X\subseteq G$ implies $\nucl(X)\subseteq \nucl(G)$, thus it suffices to study $\nucl(G)$.

$G$ acts by cellular isometries on Davis' CAT(0) realization $X_+$ of $\Delta_+$ (see \cite{Davis98}, also \cite{Caprace:2009}*{Section 2.1} and \cite{Abramenko/Brown:2008}*{Chapter 12}). For $g\in G$, denote by $|g|$ be the minimal displacement of $g$.
By \cite{Bridson:1999}*{Theorem A}, $g$ is semisimple, i.e., its minimal displacement is attained on $X_+$. By \cite{Bridson/Haefliger:1999}*{Theorem II.6.8} this implies $|g^n|= n|g|$ for $n\in\NN$.

For all $g\in \nucl(G)$ and all $n\in\NN$ there is $g_n\in G$ with $g_n^{2^n}=g$.
Hence $|g| = |g_n^{2^n}|= 2^n \cdot |g_n|$ and thus $\lim_{n\to\infty} |g_n|=0$.
But by the Proposition in \cite{Bridson:1999}, the set $\{|g| \mid g\in G\}\subseteq[0,\infty)$ is discrete.
Therefore we must have $|g|=|g_n|=0$, i.e., $g$ fixes a point in $X_+$. But that implies that $g$ stabilizes a spherical residue in $\Delta_+$. By a symmetric argument, $g$ also fixes a spherical residue in $\Delta_-$. Hence $g$ is bounded.
\end{proof}

\begin{lemma}
$A = T\cap \tau(G)$.
\end{lemma}

\begin{proof}
The inclusion $A = \tau(T) \subseteq T\cap\tau(G)$ is obvious.
Suppose now we have $g\in G$ with $\tau(g)\in T$. By the Iwasawa decomposition,
$g = bk = utk$ for some $b=ut\in B_+$, $u\in U_+$, $t\in T$, $k\in K$.
Hence $\tau(g) = u\tau(t)\theta(u)^{-1} \in U_+ A U_-$. But by the refined Birkhoff decomposition
(see \cite[Proposition~3.3(a), p.~181]{KacPeterson85c}, also \cite[Theorem~5.2.3(g)]{Kumar02}), every element of $G$ can be uniquely written as $u_+ t' u_-$ with $u_\pm \in U_\pm$ and $t'\in T$. Hence we must have $u=1$ and $\tau(g)=\tau(t)\in A$.
\end{proof}

\begin{thm} \label{nucl(tau(G)) = diagonalizable part of tau(G)}
The set
 $\nucl(\tau(G))$ equals the set of diagonalizable elements in $\tau(G)$, which in turn is the set of $K$-conjugates of $A$.
\end{thm}

\begin{proof}
By the preceding proposition, any $g\in\nucl(\tau(G))$ is bounded. By \cref{lem:sym-ss} this implies that $g$ is diagonalizable.

Suppose $g\in\tau(G)$ is diagonalizable, then it fixes some chamber $c\in\Delta_+$. But then, since $\theta(g)=g^{-1}$, it also stabilizes the chamber $\theta(c)\in\Delta_-$ opposite $c$. Since $K$ acts transitively on the pairs $(c,\theta(c))$, and since $\tau(G)$ invariant under $K$-conjugation, this implies that $g$ is contained in 
\[\bigcup_{k\in K} T^k \cap \tau(G) = \bigcup_{k\in K} (T\cap \tau(G))^k = \bigcup_{k\in K} A^k
= \bigcup_{k\in K} \nucl(A^k)\subseteq\nucl(\tau(G))
.
\qedhere
\]
\end{proof}

\section{Non-existence of polar and Cartan decompositions}

In this section we apply the results from the previous section to show that if $G$ is of non-spherical type (i.e. its Weyl group $W$ is infinite), then $G$ cannot admit a Polar or Cartan decomposition as defined in \cref{def decomps}.
Indeed, by \cref{G=tau(G)K iff nucl(tau(G))=tau(G)} we have $G=\tau(G)K$ if and only if $\nucl(\tau(G))=\tau(G)$ holds. But by \cref{nucl(tau(G)) = diagonalizable part of tau(G)}, $\nucl(\tau(G))$ consists of precisely the diagonalizable elements of $\tau(G)$. We will thus establish that $\tau(G)$ contains elements which are not diagonalizable whenever $W$ is infinite.

To illustrate why this is so, we first consider as an example a Kac-Moody group of affine, non-spherical type $\tilde{A}_n$.

\begin{example}\label{Hole}
Let $n\geq 1$ and consider the affine example $G:=\SL_{n+1}(\FF[t,t^{-1}])$ of type $\tilde{A}_n$ with the Cartan--Chevalley involution $\theta(x):=((x^{-1})^T)^\sigma$, where $\sigma$ is the $\FF$-linear ring automorphism of $\FF[t,t^{-1}]$ which interchanges $t$ and $t^{-1}$. Then let
\begin{align*}
u :=  \left(\begin{smallmatrix} 1 & 1+t \\ 0 & 1 \\ && \ddots \\ &&& 1 \end{smallmatrix}\right ) \in B_+,\qquad
v:=\tau(u) &=  u\theta(u)^{-1} =
\left(\begin{smallmatrix} 1 & 1+t \\ 0 & 1 \\ && \ddots \\ &&& 1\end{smallmatrix}\right )\cdot
\left(\begin{smallmatrix} 1 & 0 \\ 1+t^{-1} & 1 \\ && \ddots \\ &&& 1\end{smallmatrix}\right ) \\
&=
\left(\begin{smallmatrix} 1 + (1+t)(1+t^{-1}) & 1+t \\ 1+t^{-1} & 1
 \\ && \ddots \\ &&& 1\end{smallmatrix}\right )
\end{align*}
and the characteristic polynomial of $v$ is
\begin{align*}
c_\lambda(v)
&= \left( (\lambda - (1 + (1+t)(1+t^{-1}))(\lambda-1) - (1+t)(1+t^{-1}) \right) \cdot (\lambda-1)^{n-1} \\
&= \left(\lambda^2 - ( t + 4 +  t^{-1}) \lambda + 1 \right) \cdot (\lambda-1)^{n-1}.
\end{align*}
However, the polynomial $c_\lambda(v)$ does not split into linear factors over $\FF[t,t^{-1}]$, whence $v$ is not conjugate within $G$ to an element of the torus $T$, which consists of diagonal matrices with entries from $\FF$.
\end{example}

This failure to diagonalize, which can essentially be reduced to considering the Moufang tree case, i.e., type $\tilde{A}_1$, is at the heart of the general case.
While can ``fix'' this failure to split in this case by going to a suitable completion of $\FF[t,t^{-1}]$ resp. of $G$, doing so is somewhat arbitrary: There are in general multiple ways to form a completion, with different algebraic and geometric properties; moreover, we typically loose the twin building structure in the process.

\begin{lemma} \label{non-spherical => tau(G) non-diag}
Suppose $|W|=\infty$. Then $\tau(G)$ contains non-diagonalizable elements.
\end{lemma}

\begin{proof}
If the Weyl group $W$ is infinite, then by \cite{Speyer:2009}, there exists $w\in W$ such that $\ell(w^n)=n\ell(w)$ for all $n\in\NN$; such an element $w$ is called a \Defn{straight element}.
Let $c_+\in\Delta_+$ be the chambers whose stabilizer is $B_+$. It is well-known that $B_+$ acts transitively on the set $X$ of chambers in $\Delta_-$ opposite $c_+$. In an infinite building, $X$ is a connected and thick chamber system. We also have $c_-=\theta(c_+)\in X$. Therefore, there is $g\in B_+$ such that $d(c_-, g.c_-)=w$ holds. 

Now $\tau(g) c_-=g\theta(g^{-1})c_-=gc_-$ since $\theta(g^{-1})\in\theta(B_+)=B_-$.
Thus $d(c_-, \tau(g).c_-)=w$ holds. Let $\gamma$ be a minimal gallery from $c_-$ to $\tau(g).c_-$. Then, since $w$ is straight, the concatenation of $\gamma$, $\tau(g).\gamma$, $\tau(g)^2.\gamma$, \dots is still a minimal gallery.
It follows that $d(c_-,\tau(g)^nc_-)=w^n$, so $\tau(g)$ has an unbounded orbit on $G/B_-$.
But by \cref{lem:sym-ss} this means it cannot be diagonalizable.
\end{proof}

\begin{prop} \label{tauG-diag-then-spherical}
Suppose $|W|=\infty$. Then
\begin{enumerate}
\item  $\tau(G)\neq\nucl(\tau(G))$.
\item $G$ does \emph{not} admit a polar decomposition.
\item $G$ does \emph{not} admit a Cartan decomposition.
\end{enumerate}
\end{prop}

\begin{proof}
\begin{enumerate}
\item Follows from \cref{nucl(tau(G)) = diagonalizable part of tau(G)} combined with \cref{non-spherical => tau(G) non-diag}.
\item Follows from (a) and \cref{G=tauK iff tautau=tau}.

\item
Supposed we had $G=KAK$. Then $\tau(G)=\tau(KAK)
=\bigcup_{k\in K} k^{-1}\tau(A)k
\subseteq \bigcup_{k\in K} A^k$.
The elements of $A$ and also $A^k$ are diagonalizable, so all elements of $\tau(G)$ are diagonalizable.
Contradiction to \cref{non-spherical => tau(G) non-diag}.
\qedhere
\end{enumerate}
\end{proof}

The fact that there is no Cartan decomposition implies that the Kac-Moody symmetric space $G/K$ is not geodesic if $|W|=\infty$ (i.e., it contains pairs of points which are not connected by a geodesic); the following obversation implies that it is nevertheless geodesically connected (i.e., any two points can be connected by a sequence of geodesics). See also \cite{FHHK}*{Theorem 1.8}.

\begin{lemma}
Suppose $|W|=\infty$. Then
\[ G=\bigcup_{n=1}^\infty (KAK)^n. \]
\end{lemma}

\begin{proof}
Recall that $G$ is generated by its fundamental rank 1 subgroups $G_\alpha=\gen{U_\alpha,U_{-\alpha}}$. For these, the Cartan decomposition $G_\alpha=K_\alpha A_\alpha K_\alpha$ holds, where $K_\alpha:=G_\alpha\cap K$ and $A_\alpha:=G_\alpha\cap A$. The claim follows, as $G = \gen{G_\alpha \mid \alpha\in\Pi} \subseteq \bigcup_{n=1}^\infty (KAK)^n \subseteq G$.
\end{proof}

\begin{question}
Does $G=(KAK)^N$ hold for some $N\in\NN$? If so, can we bound $N$? 
\end{question}

Clearly, $N\geq 2$, but what about upper bounds? I suspect that no such $N$ exists.

In closing, we mention this ``Kostant-type'' decomposition. Geometrically, it implies that any two points in the Kac-Moody symmetric space $G/K$ are ``connected'' by a globally bounded number of horospheres.

\begin{lemma}
There is $N\in\NN$ such that $G=(KUK)^N$.
\end{lemma}

\begin{proof}
Similar to the previous proof, the fundamental rank 1 subgroups satisfy $G_\alpha = K_\alpha U_\alpha K_\alpha$.
In particular, $A_\alpha\in KUK$. Let $n$ be the rank of $G$, and $\Pi=\{\alpha_1,\dots,\alpha_n\}$ the set of fundamental roots. Then $A=A_{\alpha_1}\cdots A_{\alpha_n} \in (KUK)^n$. But then $G=UAK=AUK\subseteq (KUK)^{n+1}\subseteq G$.
\end{proof}

\begin{remark}
In the proof above, we chose $N:=n+1$, where $n$ is the rank of $G$. But we can do better: Call a \Defn{spherical covering} of the Dynkin diagram of $G$ any partition $\PPP$ of its vertices $\{1,\dots,n\}$ such that all $P\in\PPP$ corresponds to a spherical subdiagram. Clearly $\{ \{1\},\dots,\{n\}\}$ always is a spherical covering. Hence $r\leq n$ holds. A straight forward adaption of the preceding proof implies $G=(KUK)^{r+1}$.

This is still not optimal: If $G$ is of spherical type, then this gives $r=1$ and $N=2$, even though $G=KUK$ holds, i.e., one can take $N=1$.
\end{remark}

\begin{question}
Does $G=KUK$ hold when $G$ is not of spherical type?
\end{question}

\end{document}